\documentclass[10pt,a4paper]{amsart}
\usepackage[english]{babel}
\usepackage[utf8]{inputenc}
\usepackage[a4paper]{geometry}
\usepackage{amssymb,amsmath,amsthm}
\usepackage{graphicx,wrapfig,fancybox,xcolor}
\usepackage{array}
\usepackage{caption}
\usepackage[font=small]{subcaption}
\usepackage{multicol}
\usepackage{multirow}
\usepackage{tikz}
\usepackage{tikzscale}

\newtheorem{theorem}{Theorem}[section]
\newtheorem{proposition}{Proposition}[section]
\newtheorem{lemma}{Lemma}[section]
\newtheorem{question}{Question}[section]
\newtheorem{corollary}{Corollary}[section]
\newtheorem{conjecture}{Conjecture}[section]
\newtheorem{definition}{Definition}[section]
\newtheorem{remark}{Remark}[section]
\newtheorem{example}{Example}[section]

\newtheorem{problem}{Problem}[section]

\graphicspath{{images/}}

\begin{document}

\title{On Covering Simplices by Dilations in Dimensions 3 and 4}

\date{} 

\author{Lei Song, Huanqi Wen and Zhixian Zhu}

\maketitle

\begin{abstract} 
\noindent We propose a conjecture regarding the integrally closedness of lattice polytopes with large lattice lengths. We demonstrate that a lattice simplex in dimension 3  (resp. 4) with lattice length of at least 2 (resp. 3 and no edge has lattice length 5) can be covered by dilated simplices of the form $sQ$, where integer $s\ge 2$ (resp. 3) and $Q$ is a lattice simplex. The covering property implies these simplices are integrally closed. As an application, we obtain a simple criterion for the projective normality of ample line bundles on 3-(resp. 4-) dimensional $\mathbb{Q}$-factorial toric Fano varieties with Picard number one. Along the way, we discover certain unexpected phenomenon. 
\end{abstract}



\section{Introduction}
Let $M$ be a free $\mathbb{Z}$-module of rank $n\ge 1$ and $M_{\mathbb{R}}=M\otimes_{\mathbb{Z}} \mathbb{R}$ be the extension of the coefficients to the real numbers. A polytope $P$ in $M_{\mathbb{R}}$ is defined as the convex hull of a finite subset $\{u_0,u_1,\cdots,u_m\}\subset M_{\mathbb{R}}$. We refer to the polytope $P=\text{Conv}(u_0,\cdots, u_m)$ as a lattice polytope if $u_i\in M$ for all $i$. In this paper, we assume all polytopes are of full dimension. 

Let $P,\ P_1,\ P_2$ be polytopes in $M_{\mathbb{R}}$. The Minkowski sum $P_1+P_2:=\{u_1+u_2\: | \: u_1\in P_1, u_2\in P_2\}$ and the dilation by scalar $rP:=\{ru \: | \: u\in P\}$ for a positive real number $r$. When $r$ is a natural number, one can easily verify that the dilation $rP$ coincides with the $r$-fold Minkowski sum of $P$, i. e., $rP=\{u_1+\cdots+u_r \: | \: u_1,\cdots,u_r\in P\}$.

It is a fundamental problem that 
\begin{problem}
\label{prob:rP+P}
    For which lattice polytopes $P$ do the equality 
    \begin{equation}\label{eq:rP+P=(r+1)P}
        (M\cap P) + (M\cap rP) = M\cap (r+1)P
    \end{equation}
hold for all $r\in \mathbb{Z}_{>0}$?
\end{problem}

A lattice polytope that satisfies equality (\ref{eq:rP+P=(r+1)P}) for all $r\in \mathbb{Z}_{>0}$ is called to be \textit{integrally closed}. T. Oda asked in \cite{Oda08} whether any smooth lattice polytope $P$ is integrally closed. We say a polytope $P$ is smooth if, for every vertex, the primitive vectors on the edges form a basis of the lattice. In the language of algebraic geometry, Oda's question is to ask whether any ample line bundle on a smooth projective toric variety is projective normal. This problem remains widely open, and the reader is referred to \cite{miniworkshop08} for a summary of known results in this area.

To see the stated equivalence, let us recall some basics in toric geometry (cf.~\cite{Fulton93}). Consider an algebraic torus $T=\text{Spec}\mathbb{C}[M]$ of dimension $n$, where $M$ can be regarded as the character group of $T$, that is $M = \text{Hom}(T,\mathbb{C}^*)$. Given a pair $(X,L)$ with $X$ being a projective toric variety of dimension $n$ and $L$ a $T$-invariant ample line bundle on $X$, there exists an associated lattice polytope $P=P_L$ in $M_{\mathbb{R}}$. The dilation $rP$ then corresponds to the $r$-fold tensor product $L^{\otimes r}$. For $u\in M$, let $\chi^u$ denote the corresponding regular function on $T$, also viewed as a rational function on $X$. Then we have an isomorphism
\begin{equation}
    \label{eq:isomorphism}
    H^0(X,L)\cong \bigoplus_{u\in P\cap M} \mathbb{C}\cdot \chi^u.
\end{equation}
Moreover the multiplication map
\begin{equation}
    \label{eq:multiplication map}
    H^0 (X,L^{\otimes r})\otimes H^0 (X,L) \to H^0(X,L^{\otimes(r+1)})
\end{equation}
sends $\chi^{u_1}\otimes \chi^{u_2}$ for $u_1\in rP\cap M$ and $u_2\in P\cap M$ to $\chi^{u_1+u_2}$ via the isomorphism (\ref{eq:isomorphism}). Therefore the equality $(rP\cap M)+(P\cap M)=(r+1)P \cap M$ amounts to the surjectivity of the multiplication map (\ref{eq:multiplication map}). Finally recall that a base point free line bundle $L$ is said to be \textit{ projective normal} if the multiplication map (\ref{eq:multiplication map}) surjects for all $r\in \mathbb{Z}_{>0}$. It is worth mentioning projective normality of general smooth surfaces is already challenging, see \cite{GP99, Song19} and references therein.

Concerning Problem \ref{prob:rP+P}, when $n=\dim P=2$, Koelman \cite{Koelman1993Generators} established the validity of  (\ref{eq:rP+P=(r+1)P}) for all lattice polytopes $P$. For general $n$, if $P$ has an unimodular triangulation or covering, then it is integrally closed; however neither condition is necessary, see e.g. \cite{Sturmfels96, BG99}. Ogata and Nakagawa \cite{ON02} proved equality (\ref{eq:rP+P=(r+1)P}) for any $P\subset M_{\mathbb{R}}$ and $r\ge n-1$. And Ogata \cite{Ogata05} proved when $P$ is a simplex, equality (\ref{eq:rP+P=(r+1)P}) holds for $r>\frac{n-1}{2}$ provided in addition $P$ is very ample. 

Before proposing our conjecture, we first give the definition of the lattice length for lattice polytopes.
\begin{definition}
 Let $P$ be a lattice polytope in $M_{\mathbb{R}}\cong \mathbb{R}^n$. Define the lattice length of each edge $e$ of $P$ as the number of lattice points on $e$ (including end vertices) minus 1. The minimum among the lattice lengths of all edges is referred to as the lattice length of $P$, and denoted by $l(P)$. 
\end{definition}

\begin{figure}[h!] 
    \begin{subfigure}{0.45\textwidth}
        \includegraphics[width=0.8\textwidth]{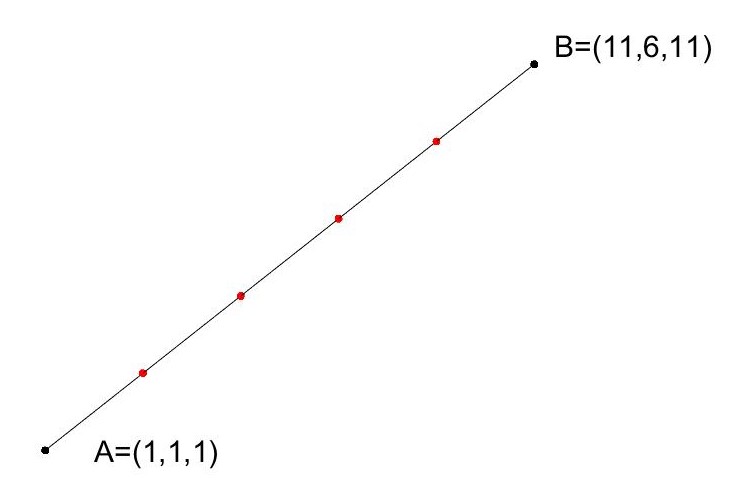}
        \caption{$l(AB)=5$}
        \label{fig:lattice length of two points}
    \end{subfigure}
    \begin{subfigure}{0.5\textwidth}
       \includegraphics[width=0.8\textwidth]{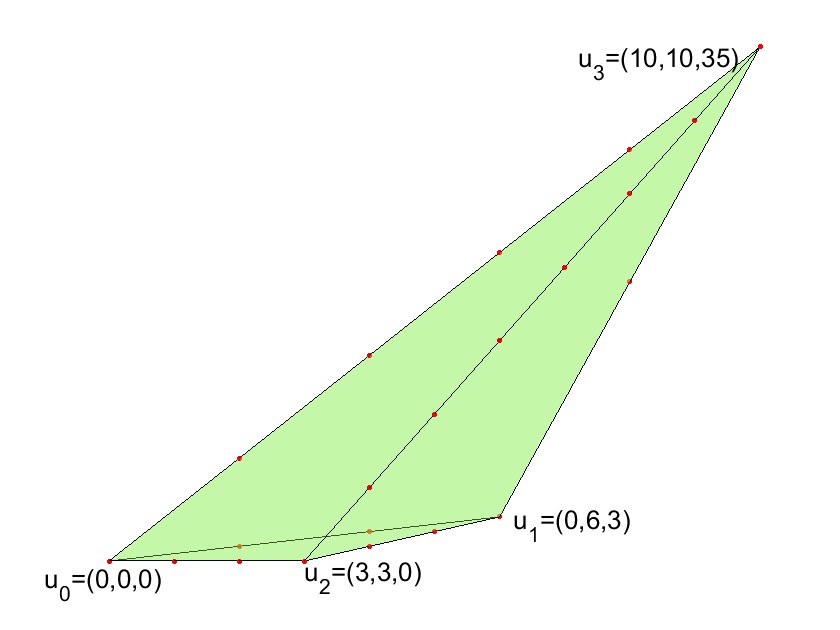}
        \caption{$l(P)=2$}
        \label{fig:lattice length of polytope}
    \end{subfigure}
    \caption{lattice length}
\end{figure}

Suppose $P=P_L$ for the pair $(X, L)$, then the lattice length of an edge $e$ can be viewed as $\deg L|_C=L\cdot C$, where $C$ is the $T$-invariant curve (hence a rational curve) associated with the 1-dimenisonal lattice polytope $e$. Thus lattice polytopes with large lattice lengths correspond to very ``positive" line bundles. Gonz{\'a}lez and the third author \cite{GZ22} showed that if $l(P)\ge n-1$, then the ample line bundle $L$ is very ample. Inspired by the Mukai conjecture (cf.~\cite{EL93}) and its variants regarding syzygies, the first and third authors made 
\begin{conjecture}[Song-Zhu]\label{S-Z Conj}
\label{conj:song&zhu}
Let $P\subset M_{\mathbb{R}}$ be a lattice polytope of dimension $n$ with $l(P)\ge n-1$, then the equality 
\[(rP\cap M) +(P\cap M)=(r+1)P\cap M\]
holds for all $r\in \mathbb{Z}_{>0}$, i.e. $P$ is integrally closed.
\end{conjecture}

By Gubeladze \cite{Gubeladze12}, Conjecture \ref{S-Z Conj} is true provided that $l(P)\ge 4n(n+1)$, and the condition can be slightly relaxed to $l(P)\ge n(n+1)$ for simplices. However it is desirable to have a condition that is linear in the dimension.

By \cite{ON02}, the conjecture is true for the dilation $P=sQ$ with $s\ge n-1$, for some lattice polytope $Q$, where $l(P)\ge n-1$ automatically. Here by $sQ$, we mean a dilation of $Q$ by $s$ up to a translation by an $u\in M$. In particular, if $P$ can be covered by $s_i Q_i$ where $s_i\ge n-1$ and $Q_i$ is a lattice simplex, then the conjecture holds for $P$. The consideration leads us to the natural question:

\begin{question}\label{quest:cover}
Given a lattice polytope $P$ with $l(P)\ge n-1$, do there exist finitely many lattice simplices $Q_i$ and integers $s_i\ge n-1$ such that $P=\bigcup s_iQ_i$?
\end{question}

If the answer to Question \ref{quest:cover} is affirmative for $P$, then Conjecture \ref{S-Z Conj} holds as well.

In this paper, we consider Question \ref{quest:cover} for lattice simplices. We shall affirmatively address this question in dimension 3, and for most cases in dimension 4. The remaining case is when $P$ has at least one edge with lattice length 5. 
\begin{theorem}\label{main thm}
A lattice simplex in dimension 3 (resp. 4) with lattice length at least 2 (resp. 3 and no edge has lattice length 5) can be covered by dilations of the form $sQ$, where integers $s\ge 2$ (resp. 3) and $Q$ are lattice simplices.
\end{theorem}

In toric Mori theory, $\mathbb{Q}$-factorial toric Fano varieties with Picard number one play an important role, and these varieties can be obtained from lattice simplices, see \cite{Fujino03}. The following is immediate. 

\begin{corollary}
Let $X$ be a  3-(resp. 4-) dimensional $\mathbb{Q}$-factorial toric Fano variety with Picard number one and $L$ be an ample line bundle on $X$. Suppose $L\cdot C\ge 2$ (resp. 3 but not equal to 5) for any $T$-invariant curve $C$ on $X$. Then $L$ is projectively normal. \qed
\end{corollary}
\textbf{Acknowledgments}.
The authors would like to thank Lujia Wang, Guoce Xin and Hanbin Zhang for helpful discussions and conversations. During the preparation of this paper, L.S. was partially supported by NSFC grants No.~12471043 and No.~12371063, and Z.Z. was partially supported by NSFC grant No.~12101423.
\vspace{2mm}

\section{Preliminaries}
In this section, we discuss specific dilations and their translations inside a lattice simplex, and fix notations throughout. The discussion is applicable to all dimensions, however we will specialize to dimensions 3 and 4 in subsequent sections. 
\subsection{Constructing dilations}
\begin{figure}[h!] 
        \makebox[0.23\textwidth]{\scriptsize $P_{0,2}$}
        \makebox[0.23\textwidth]{\scriptsize $P_{1,2}$}
        \makebox[0.23\textwidth]{\scriptsize $P_{2,2}$}
        \makebox[0.23\textwidth]{\scriptsize $P_{3,2}$}
        \\ \includegraphics[width=0.23\textwidth]{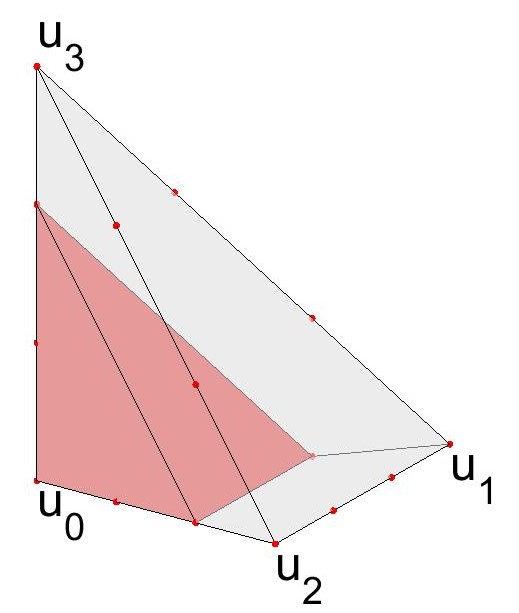}
        \includegraphics[width=0.23\textwidth]{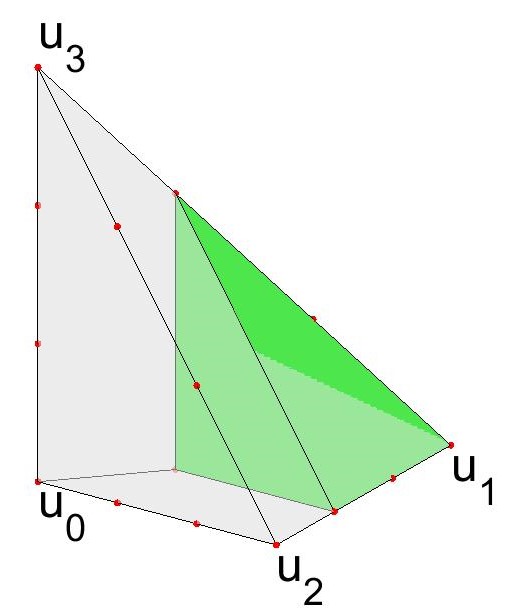}
        \includegraphics[width=0.23\textwidth]{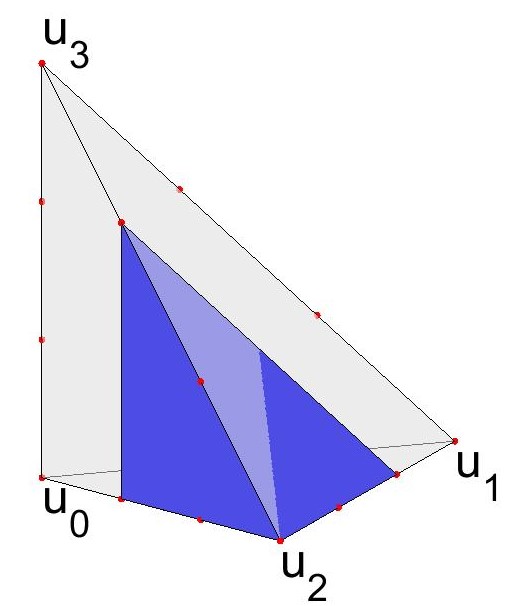}  \includegraphics[width=0.23\textwidth]{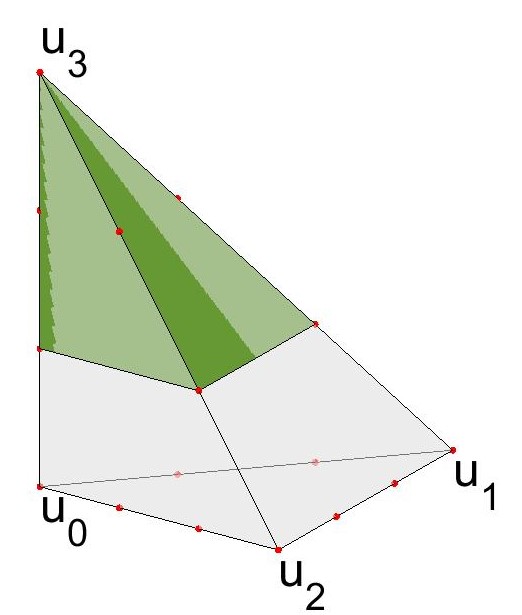}
        \caption{example: $l_{01}=l_{02}=l_{03}=l_{12}=l_{13}=l_{23}=3$}
        \label{fig:example1}
\end{figure}

\begin{figure}[h!]
        \makebox[0.23\textwidth]{\scriptsize $P_{0,3}$}
        \makebox[0.23\textwidth]{\scriptsize $P_{0,4}$}
        \makebox[0.23\textwidth]{\scriptsize $P_{0,5}$}
        \makebox[0.23\textwidth]{\scriptsize $P_{0,6}$}
        \\    \includegraphics[width=0.23\textwidth]{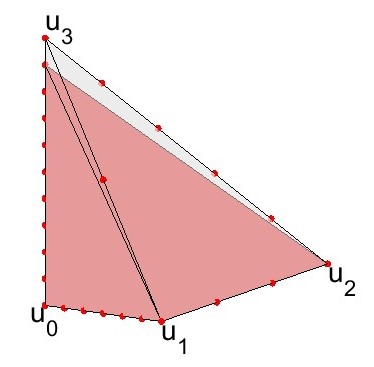}
        \includegraphics[width=0.23\textwidth]{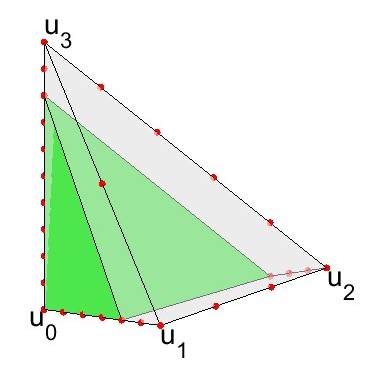}
        \includegraphics[width=0.23\textwidth]{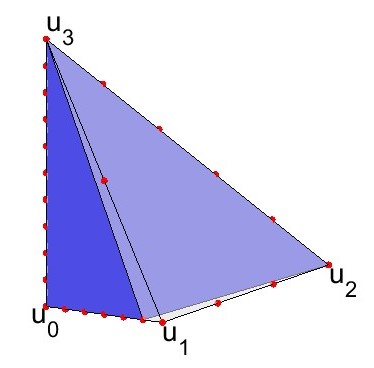}
        \includegraphics[width=0.23\textwidth]{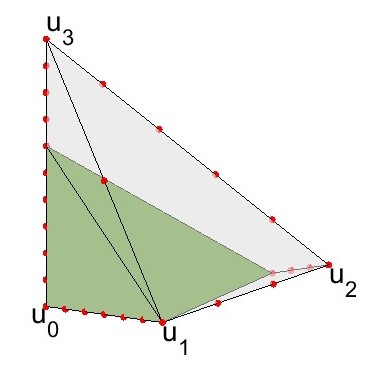}
        \caption{example: $l_{01}=6,\ l_{02}=15,\ l_{03}=10,\ l_{12}=3,\ l_{13}=2,\ l_{23}=5$}
        \label{fig:example2}
\end{figure}

Let $P=\text{Conv}(u_0,u_1,\cdots,u_n)$ be a lattice simplex of dimension $n$ and $l_{ij}$ be the lattice length of the edge from vertices $u_i$ to $u_j$. To study  Question \ref{quest:cover}, we aim to find ``maximal" lattice simplices of the form $kQ$ to cover $P$ effectively. Thus, for each $0\le i\le n$ and $2\le k \le \min_j \{l_{ij}\}$, we construct a k-dilation $P_{i,k}$ as follows: 
\begin{equation}\label{eq:V(P_i) in n-d}
    P_{i,k}=\text{Conv}(u_i,\ \frac{l_{i0}-r_{i0,k}}{l_{i0}} u_0+\frac{r_{i0,k}}{l_{i0}} u_i,\ \cdots,\ \widehat{\frac{l_{ii}-r_{ii,k}}{l_{ii}} u_i+\frac{r_{ii,k}}{l_{ii}} u_i},\ \cdots,\ \frac{l_{in}-r_{in,k}}{l_{in}} u_n+\frac{r_{in,k}}{l_{in}} u_i)
\end{equation}
where $r_{ij,k}\equiv l_{ij} (\text{mod } k)$, the least non-negative integer congruent to $l_{ij}$ modulo $k$. Some examples of $P_{i,k}$ are illustrated in Figures \ref{fig:example1} and \ref{fig:example2}. 

The following is a simple criterion to determine whether $u\in P$ is contained in $P_{i, k}$.
\begin{lemma}
    \label{le:u in P_i in n-d}
    Suppose $u\in P$, and $u$ can be uniquely expressed as $u=\sum_{i=0}^n \lambda_i u_i$, where $\lambda_i\ge 0,\ \sum_{i=0}^n \lambda_i=1$. Then  $u\notin P_{i,k}$ if and only if the inequality
\begin{equation}
        \label{eq:con of u notin p_i in n-d}
        \lambda_i < \frac{r_{i0,k}}{l_{i0}-r_{i0,k}}\lambda_0+\cdots+\widehat{\frac{r_{ii,k}}{l_{ii}-r_{ii,k}}\lambda_i}+\cdots+\frac{r_{in,k}}{l_{in}-r_{in,k}}\lambda_n
\end{equation}
holds.
\end{lemma}

\begin{proof}
   Without loss of generality, we assume $i=0$. By definition of $P_{0,k}$ in (\ref{eq:V(P_i) in n-d}),  $u\in P_{0,k}$ if and only if there exist $\mu_i$ for $0\le i\le n$ such that $\mu_i\ge 0,\ \sum_{i=0}^n \mu_i=1$ and 
    \begin{equation*}
        \sum_{i=0}^{n} \lambda_i u_i
        = \mu_0 u_0 + \sum_{i=1}^n \mu_i \left(\frac{l_{0i}-r_{0i,k}}{l_{0i}} u_i + \frac{r_{0i,k}}{l_{0i}} u_0\right).
    \end{equation*}
    Equating the corresponding coefficients of both sides of the above gives
    \begin{align}
        \lambda_0&=\mu_0 + \sum_{i=1}^n \frac{r_{0i,k}}{l_{0i}}\mu_i \label{eq2},\\
        \lambda_i&=\frac{l_{0i}-r_{0i,k}}{l_{0i}} \mu_i \quad (1\le i\le n) \label{eq3}.
    \end{align}
    Therefore, $u\in P_{0,k}$ if and only if there exist $\mu_i\ (0\le i \le n,\ \mu_i\ge 0,\ \sum_{i=0}^{n} \mu_i=1)$ such that Eqs. (\ref{eq2}) and  (\ref{eq3}) hold. We claim that the existence of such $\mu_i$ is equivalent to the following inequality
    \begin{equation}\label{eq4}
        \lambda_0 \ge \frac{r_{01,k}}{l_{01}-r_{01,k}}\lambda_1+\frac{r_{02,k}}{l_{02}-r_{02,k}}\lambda_2+\cdots+\frac{r_{0n,k}}{l_{0n}-r_{0n,k}}\lambda_n.
    \end{equation}

    Suppose the inequality (\ref{eq4}) holds, then
    \begin{align*}
        \mu_i&=\frac{l_{0i}}{l_{0i}-r_{0i,k}} \lambda_i\quad  \text{for } 1\le i\le n,\\
        \mu_0&=\lambda_0-\sum_{i=1}^n \frac{r_{0i,k}}{l_{0i}-r_{0i,k}}\lambda_i
    \end{align*}
    give the desired $\mu_i$. 
    
    Conversely, if there exist $\mu_i$ with $\mu_i\ge 0,\  \sum_{i=0}^{n} \mu_i=1$ satisfying Eqs. (\ref{eq2}) and (\ref{eq3}), then by Eq. (\ref{eq3}), we have for $1\le i\le n$,
    \[\mu_i= \frac{l_{0i}}{l_{0i}-r_{0i,k}} \lambda_i.\]
    By substituting $\mu_i\ (1\le i\le n)$ into Eq. (\ref{eq2}), we obtain
    \[\mu_0=\lambda_0-\sum_{i=1}^n \frac{r_{0i,k}}{l_{0i}-r_{0i,k}}\lambda_i.\]
    Due to $\mu_0\ge 0$, the inequality (\ref{eq4}) holds.  
\end{proof}

\subsection{Translating dilations}
In high dimensions, the constructed dilations do not suffice to cover $P$. To have more k-dilations to cover $P$, we consider translations of $P_{i,k}$. 

We fix $t_i=(t_{i0},t_{i1},\cdots,\widehat{t_{ii}},\cdots,t_{in})\in \mathbb{Z}_{\ge 0}^n$ and put $\widetilde{u_{ij}}=\frac{u_j-u_i}{l_{ij}}$. Let 
\begin{equation}\label{eq:translation of P_i,k}
P_{i,k,t_i}=P_{i,k}+\sum_{\substack{j=0\\j\ne i}}^{n} t_{ij} \widetilde{u_{ij}}.
\end{equation}

When $i=0$, simply write $t_{0}=(t_{01},\cdots,t_{0n})$ as $t=(t_1,\cdots,t_n)$, $\widetilde{u_{0j}}$ as $\widetilde{u_j}$ and $P_{0,k,t_0}$ as $P_{0,k,t}$.

\begin{figure}[h]
    \centering
    \includegraphics[width=0.95\textwidth]{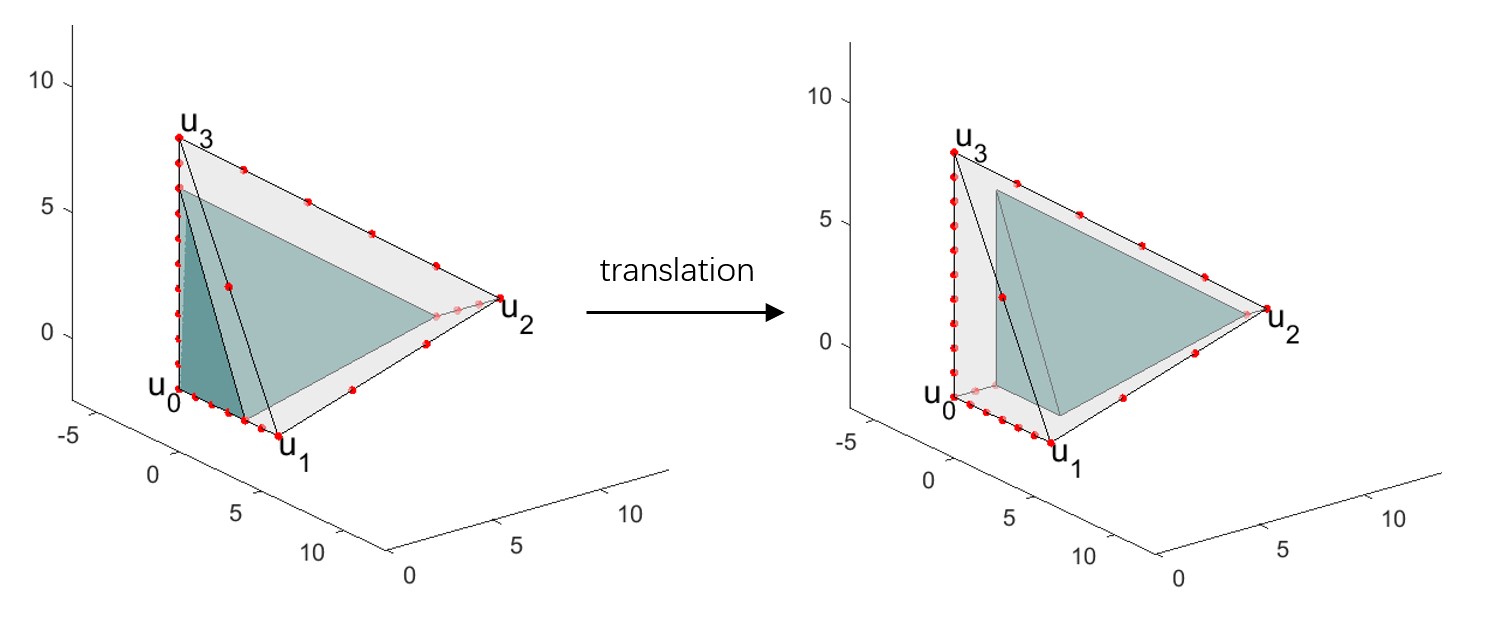}
    \caption{translation with $t=(0,2,0)$}
    \label{fig:The translation of P_{ik}}
\end{figure}

The following lemmas characterize whether the translation is still within $P$ and whether a given point is contained in the translation, respectively. 
\begin{lemma}
\label{le:P_i,k,t contained in P}
    The following statements are equivalent:
    \begin{enumerate}
        \item[(i)] $P_{0,k,t}\subseteq P$.
       \item[(ii)] For all $1\le s \le n$, $\frac{l_{0s}-r_{0s,k}}{l_{0s}} u_s+\frac{r_{0s,k}}{l_{0s}} u_0 + \sum_{j=1}^{n} t_{j} \widetilde{u_{j}}\in P.$
       \item[(iii)] $\sum_{j=1}^{n} \frac{t_j}{l_{0j}}\le \min_{1\le s\le n}\{\frac{r_{0s,k}}{l_{0s}}\}.$ 
    \end{enumerate}
\end{lemma}
\begin{proof}
For each $1\le s\le n$,
\begin{equation*}
\frac{l_{0s}-r_{0s,k}}{l_{0s}} u_s+\frac{r_{0s,k}}{l_{0s}} u_0 + \sum_{j=1}^{n} t_{j} \widetilde{u_{j}}
=\left(\frac{r_{0s,k}}{l_{0s}}-\sum_{j=1}^{n} \frac{t_j}{l_{0j}}\right) u_0+\sum_{j=1}^{n}\left(\frac{t_j}{l_{0j}}\right)u_j+\frac{l_{0s}-r_{0s,k}}{l_{0s}} u_s.
\end{equation*}
Thus $(ii)\Longleftrightarrow (iii)$. $P_{0, s, t}\subseteq P$ if and only if all its vertices are in $P$, so $(i)\Longrightarrow (ii)$. $(iii)\Longrightarrow (i)$: Note that $\sum^n_{j=1} \frac{t_j}{l_{0j}}\le 1$, which amounts to $u_0+\sum_{j=1}^{n} t_{j} \widetilde{u_{j}}\in P$.  Combining with (ii), we deduce $P_{0, k, t}\subseteq P$.
\end{proof}

\begin{lemma}
\label{le:u in P_i,k,t}
    Given $u=\sum_{i=0}^n \lambda_i u_i \in P$, then $u\in P_{0,k,t}$ if and only if
    \begin{enumerate}
        \item[(i)] For all $1\le i \le n$, 
        $\frac{t_i}{l_{0i}}\le \lambda_i$, and 
        \item[(ii)] $\frac{r_{01,k}}{l_{01}-r_{01,k}}\lambda_1+\cdots+\frac{r_{0n,k}}{l_{0n}-r_{0n,k}}\lambda_n \le \lambda_0 + \sum_{i=1}^n \frac{t_i}{l_{0i}-r_{0i,k}}.$\qed
    \end{enumerate}
\end{lemma}
We leave the proof to the interested reader. 

\section{Dimension three}
Keep the notations above.
\begin{proposition}
    \label{prop:n=3 cover}
    If $P=\text{Conv}(u_0,u_1,u_2,u_3)$ is a lattice simplex in dimension three with $l(P) \ge 2$, then there exist lattice simplices $P_{i,k}$ with $(0\le i\le 3,k\ge 2)$ such that $P=\mathop{\bigcup}\limits_{i,k} P_{i,k}$, where $P_{i,k}=kQ_{i,k}$ for some lattice polytope $Q_{i,k}$.
\end{proposition}

The polytopes $P_{i,k}$ in the proof are given by (\ref{eq:V(P_i) in n-d}). Hence it is clear that $P_{i,k}$ is a $k$-dialation of $Q_{i,k}$. To prove Proposition \ref{prop:n=3 cover}, we initially take $k=2$ in $P_{i,k}$ for all $0\le i\le 3$. Subsequently, we will assess whether $P$ can be covered by the resulting four simplices, that is, $P=\bigcup_{i=0}^3 P_{i,2}$.

According to Lemma \ref{le:u in P_i in n-d}, the following statements are equivalent:
\begin{enumerate}
    \item[(i)] $P=\mathop{\bigcup}\limits_{i=0}^3 P_{i,2}$.
    \item[(ii)] The system of inequalities for $\lambda_i$
    \begin{equation}
    \label{eq:con of u notin any P_i in 3-d}
	\left\{    
		\begin{aligned}
        &\lambda_0-\frac{r_{01,2}}{l_{01}-r_{01,2}}\lambda_1-\frac{r_{02,2}}{l_{02}-r_{02,2}}\lambda_2-\frac{r_{03,2}}{l_{03}-r_{03,2}}\lambda_3 < 0 \\
		&\lambda_1-\frac{r_{10,2}}{l_{10}-r_{10,2}}\lambda_0-\frac{r_{12,2}}{l_{12}-r_{12,2}}\lambda_2-\frac{r_{13,2}}{l_{13}-r_{13,2}}\lambda_3 < 0\\
		&\lambda_2-\frac{r_{20,2}}{l_{20}-r_{20,2}}\lambda_0-\frac{r_{21,2}}{l_{21}-r_{21,2}}\lambda_1-\frac{r_{23,2}}{l_{23}-r_{23,2}}\lambda_3 < 0\\
		&\lambda_3-\frac{r_{30,2}}{l_{30}-r_{30,2}}\lambda_0-\frac{r_{31,2}}{l_{31}-r_{31,2}}\lambda_1-\frac{r_{32,2}}{l_{32}-r_{32,2}}\lambda_2 < 0
		\end{aligned} 
	\right. 
    \end{equation}
    has no solution.
\end{enumerate}

Adding up the inequalities in (\ref{eq:con of u notin any P_i in 3-d}) yields
\begin{equation*}
    \begin{aligned}
       & \left(1-\frac{r_{10,2}}{l_{10}-r_{10,2}}-\frac{r_{20,2}}{l_{20}-r_{20,2}}-\frac{r_{30,2}}{l_{30}-r_{30,2}}\right)\lambda_0 
       +\left(1-\frac{r_{01,2}}{l_{01}-r_{01,2}}-\frac{r_{21,2}}{l_{21}-r_{21,2}}-\frac{r_{31,2}}{l_{31}-r_{31,2}}\right)\lambda_1\\
       +&\left(1-\frac{r_{02,2}}{l_{02}-r_{02,2}}-\frac{r_{12,2}}{l_{12}-r_{12,2}}-\frac{r_{32,2}}{l_{32}-r_{32,2}}\right)\lambda_2 
       +\left(1-\frac{r_{03,2}}{l_{03}-r_{03,2}}-\frac{r_{13,2}}{l_{13}-r_{13,2}}-\frac{r_{23,2}}{l_{23}-r_{23,2}}\right)\lambda_3<0.
    \end{aligned}
\end{equation*}
That is
\begin{equation}
    \label{con1:u notin any P_i in 3-d}
    \sum_{i=0}^{3} A_i \lambda_i<0,
\end{equation}
where $A_i=1-\frac{r_{0i,2}}{l_{0i}-r_{0i,2}}-\cdots-\widehat{\frac{r_{ii,2}}{l_{ii}-r_{ii,2}}}-\cdots-\frac{r_{3i,2}}{l_{3i}-r_{3i,2}}$.

Since no solution to (\ref{con1:u notin any P_i in 3-d}) implies no solution to (\ref{eq:con of u notin any P_i in 3-d}), we begin with the simpler (\ref{con1:u notin any P_i in 3-d}) to ascertain whether $P=\cup_{i=0}^3 P_{i,2}$. 

Observe that if all $A_i\ge 0$, then the inequality (\ref{con1:u notin any P_i in 3-d}) has no solution, as $\lambda_i\ge 0$ for all $i$. Thus we have
\begin{proposition}
    \label{prop:A_i>=0}
    If $A_i\ge 0$ for all $0\le i\le 3$, then $P=\mathop{\bigcup}\limits_{i=0}^3 P_{i,2}$, where $P_{i,2}$ are given by (\ref{eq:V(P_i) in n-d}).\qed
\end{proposition}

To prove Proposition \ref{prop:n=3 cover}, it remains to consider the case when some $A_i<0$. Without loss of generality, we assume that $A_0<0$. According to Table \ref{tab:n=3,mod 2}, it occurs that  at least two of $l_{10},\ l_{20}$ and $l_{30}$ must be 3, and all values are odd. 

\begin{table}[h] 
    \centering
    \caption{$\frac{r}{l-r}$ varies with $l$, $r\equiv l\: (\text{mod } 2) $}
    \begin{tabular}{c|*{11}{c}}
        \hline
          $l$           & 2 & 3          & 4 &5      &6     &7   &$\cdots$ & $2m$  &$2m+1$ &$\cdots$\\
        \hline
        $\frac{r}{l-r}$ & 0 & $\frac{1}{2}$ & 0 &$\frac{1}{4}$ & 0 & $\frac{1}{6}$ & $\cdots$ & 0 & $\frac{1}{2m}$ & $\cdots$\\
        \hline
    \end{tabular}
    \label{tab:n=3,mod 2}
\end{table}


\begin{proposition}
\label{prop:special case}
    Suppose $A_0<0$, so we assume $l_{10}=l_{20}=3$ and $l_{30}$ is odd.  Let $S$ be the set of indices $i$ where $l_{3i}\ne 2$. Then $P=(\mathop{\bigcup}\limits_{i\in S} P_{i,3})\cup P_{3,2}$.
\end{proposition}

\begin{proof}
    The proof follows a similar manner as shown previously, with the major change being the introduction of  $P_{i,3}$. Given that $l_{10}=l_{20}=3$ and $l_{30}$ is odd, we deduce that $3| l_{12}$ and that $0\in S\subseteq \{0,1,2\}$.
   
    Suppose $u=\sum^3_{i=0} \lambda_iu_i\in P\backslash\{(\mathop{\bigcup}\limits_{i\in S} P_{i,3})\cup P_{3,2}\}$. Since $u\notin P_{3,2}$, by Lemma \ref{le:u in P_i in n-d},
    \begin{equation}
    \label{eq1 in proof of special case}
        \lambda_3<\frac{r_{30,2}}{l_{30}-r_{30,2}}\lambda_0+\frac{r_{31,2}}{l_{31}-r_{31,2}}\lambda_1+\frac{r_{32,2}}{l_{32}-r_{32,2}}\lambda_2=\sum_{i\in S} \frac{r_{3i,2}}{l_{3i}-r_{3i,2}}\lambda_i.
    \end{equation}
    On the other hand, for each $i\in S$, $u\in P\setminus P_{i,3}$, then by Lemma \ref{le:u in P_i in n-d} again
    \begin{equation}
    \label{eq2 in proof of special case}
       \lambda_i<\frac{r_{i0,3}}{l_{i0}-r_{i0,3}}\lambda_0+\cdots + \widehat{\frac{r_{ii,3}}{l_{ii}-r_{ii,3}}\lambda_i} +\cdots+\frac{r_{i3,3}}{l_{i3}-r_{i3,3}}\lambda_3=\frac{r_{i3,3}}{l_{i3}-r_{i3,3}}\lambda_3.
    \end{equation}
    Substitute (\ref{eq2 in proof of special case}) into (\ref{eq1 in proof of special case}) for each $i\in S$, and we obtain
    \begin{equation}
    \label{eq3 in proof of special case}
        \lambda_3<\sum_{i\in S} \left(\frac{r_{3i,2}}{l_{3i}-r_{3i,2}}\right) \left(\frac{r_{i3,3}}{l_{i3}-r_{i3,3}}\right)\lambda_3.
    \end{equation}
    Since $\# S\le 3$, 
    \[\sum_{i\in S} \left(\frac{r_{3i,2}}{l_{3i}-r_{3i,2}}\right) \left(\frac{r_{i3,3}}{l_{i3}-r_{i3,3}}\right)\le 3\times \frac{1}{2}\times \frac{2}{3}=1,\]
   implying that (\ref{eq3 in proof of special case}) is impossible. Therefore $P=(\mathop{\bigcup}\limits_{i\in S} P_{i,3})\cup P_{3,2}$.
\end{proof}

To summarize, the proof of Proposition \ref{prop:n=3 cover} is divided into two steps. In the first step, as indicated in Proposition \ref{prop:A_i>=0}, we show that if all $A_i\ge 0$, then $P$ can be covered by four simplices of the form $2Q$. The majority of lattice simplices $P$ with large lattice lengths satisfy the condition, except in certain special cases. In the next step, we establish that $P=\bigcup_{i,k} P_{i,k}$ still holds for these special cases by introducing new dilations $P_i=3Q_i$, as presented in Proposition \ref{prop:special case}.

\begin{remark}
It might be interesting to notice that 4 dilated simplices are adequate to cover the polytope in all cases. 
\end{remark}

\section{Dimension four}
In dimension four, more dilated simplices are needed in order to fully cover a lattice simplex $P$ with $l(P)\ge 3$. However in this section we will show that the situation remains manageable as long as no edge of $P$ has lattice length $5$, by leveraging the methodology developed for dimension 3.  
\begin{proposition}
    Let $P=\text{Conv}(u_0,u_1,u_2,u_3, u_4)\subseteq M_{\mathbb{R}}\cong \mathbb{R}^4$ be a lattice simplex with  $l(P)\ge 3$. Suppose no edge of $P$ has lattice length 5, then there exist some $P_{i,k,t_i} \ (k\ge 3)$ such that $P=\mathop{\bigcup}\limits_{i,k,t_i} P_{i,k,t_i}$, where $P_{i,k,t_i}$ are given by formula (\ref{eq:translation of P_i,k}). Note that $P_{i,k}$ can be viewed as $P_{i,k,0}$.
\end{proposition}
\begin{proof}
    First, we set the modulus $k=3$ and $t_i=0$ in formula (\ref{eq:translation of P_i,k}) for all $0\le i\le 4$, Using the same argument as in the proof for dimension 3, we claim that if $u\in P\setminus\cup_{i=0}^4 P_{i,3}$, then 
    \begin{equation}
    \label{eq:con of u notin any P_i in 4-d}
    \sum_{i=0}^{4} A_i \lambda_i<0 
    \end{equation}
    where 
    \begin{equation*}
        A_i=1-\frac{r_{0i,3}}{l_{0i}-r_{0i,3}}-\cdots-\widehat{\frac{r_{ii,3}}{l_{ii}-r_{ii,3}}}-\cdots-\frac{r_{4i,3}}{l_{4i}-r_{4i,3}}.
    \end{equation*}
    If $A_i\ge 0$ for all $i$, then the inequality (\ref{eq:con of u notin any P_i in 4-d}) clearly has no solution for $\lambda_i\ge 0$, and hence $P=\cup_{i=0}^4 P_{i,3}$.
   
    Keeping in mind that no edge of $P$ has lattice length 5, we conclulde that $A_i\ge 0$ unless there are at least two instances of $\frac{1}{3}$ in the fractions of the form $\frac{r}{l-r}$ by looking at Table \ref{tab:n=4,mod 3}. Thus we only need to focus on the situation where there are at least two edges starting from one vertex $u_i$ with lattice length 4 or 8 such that $A_i<0$. Assuming $i=0$, we proceed to discuss by cases based on the number of edges starting from $u_0$ with lattice length 4 or 8. 
    \begin{table}[h]
    \caption{$\frac{r}{l-r}$ varies with $l$, $r\equiv l\: (\text{mod } 3)$}
    \label{tab:n=4,mod 3}
    \small
    \centering
        \begin{tabular}{c|*{20}{c}}
        \hline
          $l$           & 3 & 4    &  5       &6   &7     &8   &9   & 10   &11  &12 &13 &14 &15 &16 &17 &$\cdots$ &$3m$ &$3m+1$ & $3m+2$ 
         &$\cdots$ \\
        \hline
        $\frac{r}{l-r}$ & 0 & $\frac{1}{3}$ & $\frac{2}{3}$ &0 & $\frac{1}{6}$ &$\frac{2}{6}$ &0 & $\frac{1}{9}$ &$\frac{2}{9}$ &0 & $\frac{1}{12}$ & $\frac{2}{12}$ &0 &$\frac{1}{15}$ & $\frac{2}{15}$ &$\cdots$  &0 &$\frac{1}{3m}$ & $\frac{2}{3m}$ &$\cdots$\\
        \hline
        \end{tabular}
    \end{table}
    
        \textbf{Case A}: If the lattice length of all edges staring from $u_0$ is equal to 4 or 8, then all edges of $P$ have lattice length divisible by 4. In other words, $P=4Q$ for some lattice simplex $Q$. We are done by \cite{ON02}.
        
        \textbf{Case B}: There are exactly three edges starting from $u_0$ with lattice length 4 or 8. We assign $l_{01}=4n_1$, $l_{02}=4n_2$ and $l_{03}=4n_3$ where $n_i=1$ or $2$. We note the following facts:
        \begin{enumerate}
            \item Edges $u_1u_2,\ u_1u_3$ and $u_2u_3$ have lattice length divisible by 4.
            \item Since $A_0<0$, $l_{04}$ cannot be divisible by 3. 
            \item We can assume other edges $u_4u_i$ for $0\le i\le 3$ have lattice length not divisible by 4; otherwise, it reduces to Case A: $P=4Q$.
        \end{enumerate}
        \begin{figure}[h]
        \centering
        \includegraphics[width=0.3\textwidth]{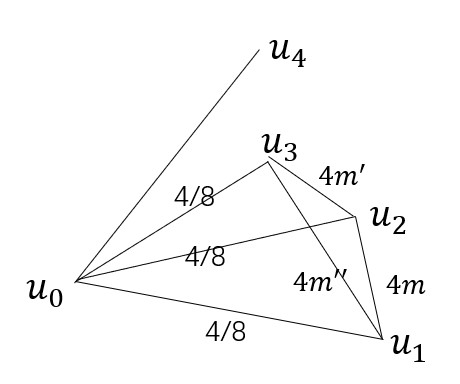}
        \caption{Case B}
        \label{fig:Case B in 4-d}
        \end{figure}
        Now let $S=\{i\ |\ l_{i4}\ne 3\}$. Since $l_{04}$ cannot be divisible by 3, $0\in S\subseteq\{0,1,2,3\}$. It follows that 
            \begin{equation*}
            \label{eq1 in Case B.2}
                P=(\cup_{i\in S} P_{i,4})\cup P_{4,3}
            \end{equation*}
        from an analogous argument as in the proof of Proposition \ref{prop:special case}.

        \textbf{Case C}: There are exactly two edges starting from $u_0$ with lattice length 4 or 8. Let the two edges be $u_0u_1$ and $u_0u_2$. Since $A_0<0$ and $l_{03},\ l_{04}$ cannot be divisible by 4, $l_{03}$ and $l_{04}$ have one value of $11$. We assume $l_{03}=11$. Then the value of $l_{04}$ have four possibilities: $7, 11, 14, 17$. We can infer some information about the lattice length of other edges: (1) $l_{12}$ is divisible by 4; (2) $l_{13}$ and $l_{23}$ are not divisible by 4 or 11; (3) $l_{14}$ and $l_{24}$ cannot be divisible by 4. Keep these in mind, as we will always use them to estimate coefficients.
        \begin{figure}[h]
        \centering
        \includegraphics[width=0.3\textwidth]{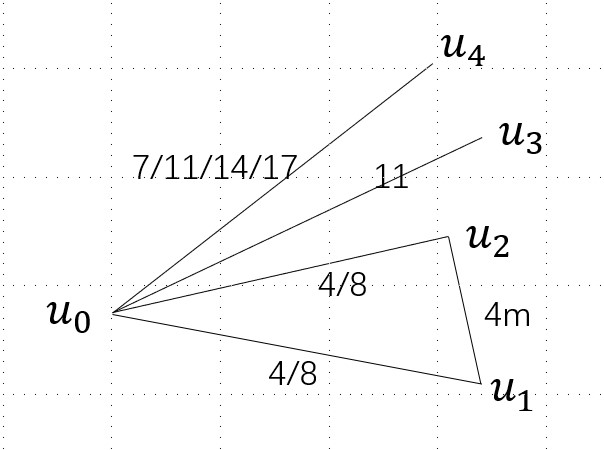}
        \caption{Case C}
        \label{fig:Case C in 4-d}
        \end{figure}

        In this case, we need to consider translations of constructed dilations to fully cover the simplex. In view of Lemma \ref{le:P_i,k,t contained in P}, to get a legitimate translation of $P_{0, 3}$, we need to find $t=(t_1, t_2, t_3, t_4)\in \mathbb{Z}_{\ge 0}^4$ such that 
        \[\sum_{j=1}^4 \frac{t_j}{l_{0j}}\le \min_{1\le s\le 4} \{\frac{r_{0s,3}}{l_{0s}}\}.\]
        It is easy to verify that $t=(t_1,t_2,t_3,t_4)=(0,0,0,1)$ and $t^{\prime}=(t_1^{\prime},t_2^{\prime},t_3^{\prime},t_4^{\prime})=(0,0,1,0)$ are two solutions, given $l_{01},l_{02}\in\{4,8\},\ l_{03}=11,\ l_{04}\in\{7,11,14,17\}$. Applying Lemma \ref{le:u in P_i,k,t} and Lemma \ref{le:u in P_i in n-d}, we obtain the following necessary and sufficient conditions for a point $u$ to lie outside the constructed dilations and the translations: 
        \begin{eqnarray}\label{eq: u notin P_0,3,t}
            u\in P\setminus P_{0,3,t} &\Longleftrightarrow&
            \lambda_4<\frac{1}{l_{04}} \text{ or } \\
            \notag
            &&\lambda_0+\frac{1}{l_{04}-r_{04,3}}<\frac{1}{3}\lambda_1+\frac{1}{3}\lambda_2+\frac{2}{9}\lambda_3+\frac{r_{04,3}}{l_{04}-r_{04,3}}\lambda_4;
        \end{eqnarray}
        \begin{eqnarray}\label{eq: u notin P_0,3,t'}
            u\in P\setminus P_{0,3,t^{\prime}}&\Longleftrightarrow&
            \lambda_3<\frac{1}{11} \text{ or }\\
            \notag
            &&\lambda_0+\frac{1}{9}<\frac{1}{3}\lambda_1+\frac{1}{3}\lambda_2+\frac{2}{9}\lambda_3+\frac{r_{04,3}}{l_{04}-r_{04,3}}\lambda_4;
        \end{eqnarray}
        \begin{equation}\label{eq:u notin P_0,4}
            u\in P\setminus P_{0,4}\Longleftrightarrow
            \lambda_0<\frac{3}{8}\lambda_3+\frac{r_{04,4}}{l_{04}-r_{04,4}}\lambda_4;
        \end{equation} 
        for $0\le i\le 4$,
        \begin{equation}\label{eq:u notin P_i,3}
            u\in P\setminus P_{i,3}\Longleftrightarrow
            \lambda_i < \frac{r_{i0,3}}{l_{i0}-r_{i0,3}}\lambda_0+\cdots+\widehat{\frac{r_{ii,3}}{l_{ii}-r_{ii,3}}\lambda_i}+\cdots+\frac{r_{i4,3}}{l_{i4}-r_{i4,3}}\lambda_4.
        \end{equation}
        Summing (\ref{eq:u notin P_i,3}) over $i$, we obtain a necessary condition 
        \begin{equation*}\label{eq:u notin all P_i,3}
            \sum_{i=0}^{4} A_i \lambda_i<0
        \end{equation*}
        for $u\in P\setminus\cup_{i=0}^4 P_{i,3}$ which is same as (\ref{eq:con of u notin any P_i in 4-d}).
        
       \textbf{Claim}: $P=\left(\cup_{i=0}^4 P_{i,3}\right)\cup P_{0,4}\cup P_{0,3,t}\cup P_{0,3,t^{\prime}}.$
        
        Suppose to the contrary that $P$ is not fully covered by these simplices. Then there exist $\lambda_i\ge 0$ with $\sum^4_{i=0}\lambda_i=1$ satisfying conditions (\ref{eq: u notin P_0,3,t}), (\ref{eq: u notin P_0,3,t'}), (\ref{eq:u notin P_0,4}) and (\ref{eq:u notin P_i,3}) simultaneously.  We will argue the impossibility of such a scenario in three cases where both conditions (\ref{eq: u notin P_0,3,t}) and (\ref{eq: u notin P_0,3,t'}) are satisfied: 
        
            (i). Suppose that $\lambda_4<\frac{1}{l_{04}}$ in (\ref{eq: u notin P_0,3,t}) and $\lambda_3<\frac{1}{11}$ in (\ref{eq: u notin P_0,3,t'}) occur. Then we have
            \begin{equation}
                \lambda_1+\lambda_2=1-\lambda_0-\lambda_3-\lambda_4>1-\frac{11}{8}\lambda_3-\left(\frac{r_{04,4}}{l_{04}-r_{04,4}}+1\right)\lambda_4>\frac{5}{8}
            \end{equation}
            where the first inequality follows from (\ref{eq:u notin P_0,4}), and the second one holds because $l_{04}\in \{7,11,14,17\}$. On the other hand, by setting $i=1,2$ in (\ref{eq:u notin P_i,3}), calculating the upper bounds of each coefficient, and summing the two resulting inequalities, we derive 
            \begin{equation*}
                \lambda_1+\lambda_2<\lambda_0+\frac{1}{2}\lambda_3+\frac{2}{3}\lambda_4<\frac78 \lambda_3+\left(\frac{r_{04,4}}{l_{04}-r_{04,4}}+\frac23\right)\lambda_4<\frac{5}{8},
            \end{equation*}
            which leads to a contradiction. Consequently, it is impossible that both $\lambda_3<\frac{1}{11}$ in (\ref{eq: u notin P_0,3,t}) and $\lambda_4<\frac{1}{l_{04}}$ in (\ref{eq: u notin P_0,3,t'}) hold simultaneously.
            
            (ii).  Suppose $\lambda_0+\frac{1}{l_{04}-r_{04,3}}<\frac{1}{3}\lambda_1+\frac{1}{3}\lambda_2+\frac{2}{9}\lambda_3+\frac{r_{04,3}}{l_{04}-r_{04,3}}\lambda_4$ in (\ref{eq: u notin P_0,3,t}) holds. Since this inequality is stronger than (\ref{eq:u notin P_i,3}) when $i=0$, we can substitute this strengthened inequality for formula (\ref{eq:u notin P_i,3}) with $i=0$ and derive an enhanced condition
            \begin{equation*}
                \frac{1}{l_{04}-r_{04,3}}+\sum_{i=0}^4 A_i\lambda_i<0.
            \end{equation*}
            However, we assert that this is impossible. In fact, we have
            \begin{align*}
                0&>\frac{1}{l_{04}-r_{04,3}}+\sum_{i=0}^4 A_i\lambda_i\\
                &>\frac{1}{l_{04}-r_{04,3}}+A_0\left(\frac{3}{8}\lambda_3+\frac{r_{04,4}}{l_{04}-r_{04,4}}\lambda_4\right)+\sum_{i=1}^4 A_i\lambda_i\\
                &=\frac{1}{l_{04}-r_{04,3}}+A_1\lambda_1+A_2\lambda_2+\left(\frac{3}{8}A_0+A_3\right)\lambda_3+\left(\frac{r_{04,4}}{l_{04}-r_{04,4}}A_0+A_4\right)\lambda_4\\
                &\ge \frac{1}{l_{04}-r_{04,3}}-\frac{1}{18}\lambda_1-\frac{1}{18}\lambda_2+\frac{5}{72}\lambda_3-\frac{1}{9}\left(\frac{r_{04,4}}{l_{04}-r_{04,4}}\right)\lambda_4\\
                &>\frac{1}{l_{04}-r_{04,3}}(1-\lambda_1-\lambda_2-\lambda_3-\lambda_4)\\
                &\ge 0.
            \end{align*}
            The second inequality holds because $A_0<0$ and (\ref{eq:u notin P_0,4}) is satisfied. The third inequality is ensured by the non-negativity of $\lambda_i$ and the lower bounds on $A_i$: $A_0\ge -\frac{1}{9}$, $A_1\ge -\frac{1}{18}$, $A_2\ge -\frac{1}{18}$, $A_3\ge \frac{1}{9}$ and $A_4\ge 0$---all are determined by estimating the upper bounds of each fraction $\frac{r}{l-r}$ in the expression of $A_i$. The subsequent inequality is guaranteed by the non-negativity of $\lambda_i$ and $\frac{1}{l_{04}-r_{04,3}}>\frac{1}{9}\cdot\frac{r_{04,4}}{l_{04}-r_{04,4}}$ when $l_{04}\in \{7,11,14,17\}$.

            (iii). Suppose $\lambda_0+\frac{1}{9}<\frac{1}{3}\lambda_1+\frac{1}{3}\lambda_2+\frac{2}{9}\lambda_3+\frac{r_{04,3}}{l_{04}-r_{04,3}}\lambda_4$ in (\ref{eq: u notin P_0,3,t'}) holds. By applying the same argument as in (ii), but replacing $\frac{1}{l_{04}-r_{04,3}}$ with $\frac{1}{9}$, we can still arrive at a contradiction.
            
  This concludes the proof of the claim, and hence the proof of the proposition.        
\end{proof}

\section{Numerical Experiment: An amusing Example}
In the preceding section, we addressed Question \ref{quest:cover} in dimension 4, but with the situation where the lattice simplex $P$ has at least one edge with lattice length 5 untouched. In this section, we illustrate through an amusing example that in such a scenario, the constructed dilations and their translations may not be adequate to cover the simplex $P$, so more dilations within $P$ need to be taken into account. 

\begin{example}
\label{example: special case with lattice length 5}
    Let $P$ be the simplex of dimension 4 with vertices 
    \[u_0=(5,0,0,0),\ u_1=(0,60,0,0),\ u_2=(0,0,0,0),u_3=(8,24,12,0),\ u_4=(33,24,72,60).\]
    The lattice lengths of edges are given by
    \[(l_{ij})_{5\times 5}=\begin{bmatrix}
         & 5 & 5 & 3 &4\\
         5 &  & 60 & 4 & 3\\
         5 & 60 &  & 4 &3 \\
         3 & 4 & 4 &  & 5\\
         4 & 3 & 3 & 5& 
    \end{bmatrix}.\]
    
    Using the previous method, we can construct five simplices $P_{i,3}\ (0\le i \le 4)$ of the form $P_{i,3}=3Q_{i,3}$. Notice that in the matrix of lattice length, every row contains at least one value of 3. According to Lemma \ref{le:P_i,k,t contained in P}, these $P_{i,3}$ cannot be translated without getting out of $P$.
    
    Then we consider whether $P_{i,3}$ can cover $P$, that is, whether the system of inequalities
    \begin{equation}
    \label{eq in example}
        \left\{\begin{aligned}
            \lambda_0<&   &\frac23 \lambda_1&+\frac23 \lambda_2&        &+\frac13 \lambda_4\\
            \lambda_1<&\frac23 \lambda_0& 
 &   &+\frac13\lambda_3&  \\
            \lambda_2<&\frac23\lambda_0& 
 &    &+\frac13\lambda_3&   \\
            \lambda_3<& 
 &\frac13\lambda_1&+\frac13\lambda_2& 
 &+\frac23\lambda_4\\
            \lambda_4<&\frac13\lambda_0& 
 &   &+\frac23\lambda_3&  
        \end{aligned}\right.
    \end{equation}
    has a solution for $\lambda_i$ such that $0\le \lambda_i\le 1$ and $\sum_{i=0}^4 \lambda_i=1$.  For instance, $u=(\lambda_0,  \cdots, \lambda_4)=(\frac{23}{100}, \frac{18}{100}, \frac{18}{100}, \frac{23}{100}, \frac{18}{100})$ is a solution, and hence $u \notin \bigcup^4_{i=0} P_{i, 3}$. 

To calculate the area not covered by $P_{i,3}\ (0\le i\le 4)$, we use the Monte-Carlo method: take a random point $(\lambda_0,\cdots, \lambda_4)$ in the standard simplex and test whether $\lambda_i,\ 0\le i\le 4$ satisfy the system (\ref{eq in example}). With the increase in the number of sampling points from $10^4$ to $10^8$, the rate of sampling points satisfying (\ref{eq in example}) stabilizes at around $1.1\%$. This indicates that the subarea of $P$ not covered by $P_{i,3},\ 0\le i\le 4$ accounts for approximately $1.1\%$.
    
\end{example}

Although $P$ in the example cannot be covered by the k-dilations $P_{i,k}$ along with their translations $P_{i,k,t}$, we still expect the general validity of Conjecture \ref{S-Z Conj}; indeed, in Example \ref{example: special case with lattice length 5}, we can identify additional k-dilations contained in $P$ to cover the area left uncovered by $P_{i,3}$. With the assistance of computer, we enumerate all lattice simplices $Q$ such that $3Q$ is contained in $P$ and identify three specific 3-dialations: $\widetilde{P_1}, \widetilde{P_2}, \widetilde{P_3}$. These, in conjunction with $P_{i, 3}\ (0\le i\le 4)$, fully cover $P$. The $\widetilde{P_i}$ are given by
\begin{eqnarray*}
\widetilde{P_1}&=&\text{Conv}( (2, 0,0,0),\ 
     (2,    42   ,  3 ,    0),\ 
     (8 ,   24  ,  12  ,   0),\ 
    (26  ,  18 ,   54   , 45),\ 
     (5   ,  0,     0    , 0)),\\
\widetilde{P_2}&=&\text{Conv}(
  (2 ,    0    , 0,     0),\ 
    ( 2,     3   ,  0,     0),\ 
    (11 ,   33  ,  21 ,   12),\ 
    (20  ,  27 ,   42  ,  33),\ 
     (5   ,  0,     0   ,  0)),\\
\widetilde{P_3}&=&\text{Conv}(
     (2    , 0,     0    , 0),\ 
     (2   ,  3 ,    0   ,  0),\ 
     (8  ,  24  ,  12  ,   0),\ 
    (14 ,   33   , 30 ,   24),\ 
     (5,     0    , 0,     0)).
\end{eqnarray*}
Similar to Lemma \ref{le:u in P_i in n-d}, a point $u=\sum_{i=0}^4 \lambda_i u_i\in P$ is not in $\widetilde{P_i}$ for any $1\le i\le 3$ if and only if $\lambda_i$ satisfy the following inequalities: 
\begin{equation}
    \label{eq: u not in widetildeP}
        \left\{\begin{aligned}     \lambda_0<&-\frac16\lambda_1&+\frac23\lambda_2&   &+\frac13\lambda_4\\
\lambda_0<&\frac23\lambda_1&+\frac23\lambda_2 &-\frac27\lambda_3&-\frac{1}{21}\lambda_4\\
\lambda_0<&\frac23\lambda_1&+\frac23\lambda_2&   &-\frac56\lambda_4
    \end{aligned}\right.
\end{equation}

Thus, $u$ is neither in $\widetilde{P_i}\ (1\le i\le 3)$ nor in $P_{i,3}\ (0\le i\le 4)$ if and only if both (\ref{eq in example}) and (\ref{eq: u not in widetildeP}) hold. To show the non-existence of such a point $u$, we add ``boundary" to the inequalities (\ref{eq in example}) and (\ref{eq: u not in widetildeP}) and introduce an objective function $\min f=\lambda_0$. This transforms the problem into a linear programming problem so that one can employ scientific software such as MATLAB to solve it. It turns out no feasible solution exists to the LP problem, so no values for $\lambda_i$ satisfy inequalities  (\ref{eq in example}) and (\ref{eq: u not in widetildeP}). This indicates that $P$ can be fully covered.

\bibliography{On_covering_simplices_by_dilations}{}
\bibliographystyle{plain}	

\bigskip
\noindent\small{\textsc{School of Mathematics, Sun Yat-sen University\\
W. 135 Xingang Rd., Guangzhou, Guangdong 510275, P.R.~China}\\
\emph{E-mail address}:  \texttt{songlei3@mail.sysu.edu.cn}

\bigskip
\noindent\small{\textsc{School of Mathematics, Sun Yat-sen University\\
W. 135 Xingang Rd., Guangzhou, Guangdong 510275, P.R.~China}\\
\emph{E-mail address}:  \texttt{wenhq7@mail2.sysu.edu.cn}

\bigskip
\noindent\small{\textsc{Academy for Multidisciplinary Studies, Capital Normal University\\
No. 105 West 3rd Ring Road, Beijing 100048, P.R.~China}\\
\emph{E-mail address}:  \texttt{zhixian@cnu.edu.cn}

\end{document}